\documentclass[12pt]{amsart}
\usepackage{graphicx}

\newtheorem{thm}{Theorem}[section]
\newtheorem{lemma}[thm]{Lemma}
\newtheorem{defin}[thm]{Definition}
\newtheorem{conjecture}[thm]{Conjecture}
\newtheorem{prop}[thm]{Proposition}
\newtheorem{rem}[thm]{Remark}
\newtheorem{cor}[thm]{Corollary}

\newcommand{\B}{{\mathbf B}}

\newcommand{\I}{{\bf e}}

\newcommand{\Hom}{{\rm Hom}}
\newcommand{\C}{\Bbb{C}}
\newcommand{\CC}{\mathbf C}

\newcommand{\R}{\Bbb{R}}
\newcommand{\Z}{\Bbb{Z}}

\newcommand{\fk}{\mathfrak{k}}
\newcommand{\G}{\mathcal{G}}

\newcommand{\fg}{\mathfrak{g}}
\newcommand{\fh}{\mathfrak{h}}

\newcommand{\im}{{\rm im}}

\renewcommand{\S}{{\mathbf S}}

\newcommand{\SL}{{\rm SL}}

\newcommand{\PGL}{{\rm PGL}}

\newcommand{\Aut}{{\rm Aut}}

\newcommand{\pr}{{\mathbf{pr}}}
\renewcommand{\H}{{\rm H}}

\newcommand{\Ad}{{\rm Ad}}

\newcommand{\HH}{{\mathsf{H}}}

\newcommand{\f}{{\mathbf f}}
\newcommand{\V}{{\mathbf V}}

\newcommand{\X}{{\mathbf X}}
\newcommand{\tSL}{\widetilde{{\rm SL}(2,\R)}}

\title{Dehn twists and invariant classes}
\commby{Daniel Ruberman}

\author{Eugene Z. Xia}
\address{Department of Mathematics, National Cheng Kung University and National Center for Theoretical Sciences, Tainan 701, Taiwan}
\email{ezxia@ncku.edu.tw}

\thanks{The author gratefully acknowledges partial support by the National Science Council, Taiwan with grants 96-2115-M-006-002 and 97-2115-M-006-001-MY3.}

\begin{document}

\begin{abstract}
A degeneration of compact K\"ahler manifolds gives rise to a monodromy action on Betti moduli space $$\HH^1(X, G) = \Hom(\pi_1(X),G)/G$$ over smooth fibres with a complex algebraic structure group $G$ being either abelian or reductive.  Assume that the singularities of the central fibre is of normal crossing.  When $G = \C$, the invariant cohomology classes arise from the global classes.  This is no longer true in general.  In this paper, we produce large families of locally invariant classes that do not arise from global ones for reductive $G$.  These examples exist even when $G$ is abelian, as long as $G$ contains multiple torsion points.  Finally, for general $G$, we make a new conjecture on local invariant classes and produce some suggestive examples.
\end{abstract}

\maketitle
\section{Introduction}
Let $X$ be a compact K\"{a}hler manifold and $G$ an abelian or reductive complex algebraic group.  The adjoint action of $G$ on the representation variety $\Hom(\pi_1(X),G)$ gives rise to the categorical quotient \cite{LM1, M1}
$$
\HH^1(X,G) = \Hom(\pi_1(X),G)/G.
$$
The functor $\HH^1(\cdot, G)$ is contravariant and equals the usual first Betti cohomological functor when $G = \C$.  Let
$$
\pr : \Hom(\pi_1(X),G) \to \HH^1(X,G)
$$
be the canonical projection and for a representation $\rho \in \Hom(\pi_1(X),G)$, let $[\rho] = \pr(\rho).$

Let $\f : \X \to D$ be a proper holomorphic map, where $\X$ is a K\"ahler manifold and $D \subset \C$ is the unit open disk.  Let $D^* = D \setminus \{0\}$.  If $\f$ has maximum rank over all points of $D^*$ and $\f^{-1}(0)$ has normal crossing singularitie(s), then $\f$ is called a {\bf degeneration of K\"ahler manifolds} (or {\bf degeneration} for short).  Fix $s \in D^*$ once and for all and let  $X = \f^{-1}(s)$.
This setting gives rise to the Picard-Lefschetz diffeomorphism \cite{C1}
\[
T: X \to X.
\]

Let $X_0 = \f^{-1}(0)$.  Then $\f$ induces a strong deformation retraction $\X \to X_0$.  This together with the inclusion $X \hookrightarrow \X$ give a map
\[
c : X \hookrightarrow \X \to X_0.
\]
Define
\begin{equation}
c^* = \HH^1(c,G) \text{ and } T^* = \HH^1(T,G).
\end{equation}\label{cT}
The local invariant cycle theorem implies that if $G = \C$, then
\begin{equation} \label{csa}
\im(c^*) = \HH^1(X,\C)^{T^*},
\end{equation}
where $\HH^1(X,\C)^{T^*}$ denotes the subset of $\HH^1(X,\C)$ fixed by $T^*$ \cite{C1, PS1, S1, S2}.
This is no longer true when $G$ is non-abelian.  Isolated examples were found for $G = \SL(n,\C)$ for various $n$ \cite{TX1}.  These examples are of the type $[\rho] \in \HH^1(X,G)$ with finite $\im(\rho)$.

For the rest of the paper, assume $G$ to be reductive; moreover, for a given degeneration $\f : \X \to D$, always assume $X = \f^{-1}(s)$ and that it is a compact Riemann surface with genus $p > 0$ and that $X_0 = \f^{-1}(0)$ is singular with normal crossing.  Let $\I$ be the identity element and $N$ the subset of torsion points of $G$.  Since $G$ is reductive, $N \neq \{\I\}$.
Let $$\V_\f =  \HH^1(X,G)^{T^*} \setminus \im(c^*).$$
We refer to $\V_\f$ as the exceptional family of $\f$.
Theorem \ref{thm:exception} shows that $\V_\f$ can be rather large for the simple reason of $N \neq \{\I\}$ and this is true even when $G$ is abelian (such as $\C^\times$).
From this perspective, one has Equation \eqref{csa} because $\C$ is unipotent with $0$ being its only torsion point.

We introduce the more restrictive notion of simple degeneration (Definition \ref{def: simple degeneration}) and construct more examples with respect to simple degenerations with Theorem \ref{thm:simple exception}.

These examples motivate us to modify the local invariant class statement to Conjecture \ref{conj:lic}.
We define the notion of pseudo-degeneration (Definition \ref{def:pseudo-degeneration}) and produce exceptional families of pseudo-degenerations.  The existence of these families and the example in Section~\ref{sec:tSL} shed light on Conjecture \ref{conj:lic}.

\vskip 0.2in
\noindent {\bf Acknowledgement:} The author is grateful to Domingo Toledo for helpful discussions, to Jiu-Kang Yu for reading previous versions of the manuscript and for comments and discussions and to the referee for careful review.

\section{Definitions and constructions}
\subsection{The fundamental groups and their representations}
We begin by recalling our standing assumptions that a degeneration $\f : \X \to D$ always has singular central fibre $X_0 = \f^{-1}(0)$ with normal crossing and, for $s \in D^*$, $X=\f^{-1}(s)$ is a compact Riemann surface of genus $p>0$ through out the rest of the paper.  Moreover,
$$\V_\f =  \HH^1(X,G)^{T^*} \setminus \im(c^*).$$

%

Let $x \in X$ and $\pi(x) = \pi_1(X,x)$ be the fundamental group of $X$ with base point $x$.  The definition of $\HH^1(X,G)$ depends on $\pi(x)$, but is independent of $x$.  It is often true in this subject that while the technical calculation in $\pi(x)$ involves $x$, the results seldom depend on $x$.
We shorten $\pi(x)$ to $\pi$ and largely ignore the issue when the base point is irrelevant.  We do not distinguish an element in $\pi$ and its representative directed loop in $X$.
All loops are assumed to be smooth.  
A loop in $X$ is {\bf simple} if it has no self-intersection; moreover, it is called {\bf simple and separating} if it is null-homologous in $\H_1(X,\Z)$; otherwise, it is called {\bf simple and non-separating}.

\vskip 0.2in
\centerline{\includegraphics[scale=0.7]{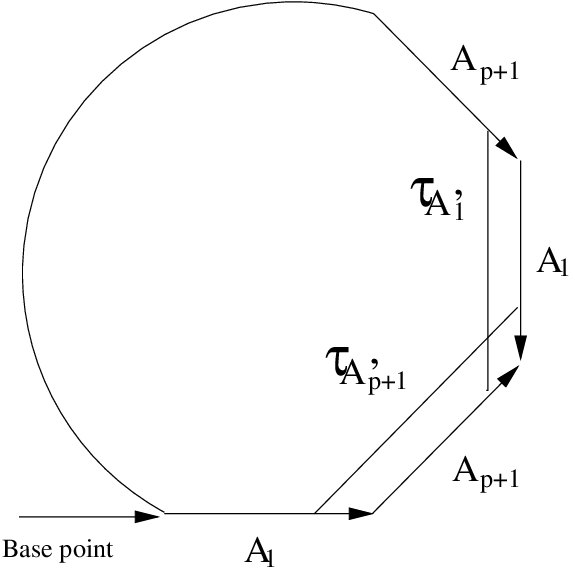}}

\centerline{F{\sc igure} 1: The generators of $\pi$ and Dehn twists along them}
\vskip 0.2in

Two concepts are ubiquitous, the adjoint action and the commutator map.
\begin{defin}
Let $\G$ be any group.  The dot notation $g.f$ means $gfg^{-1}$ where $g \in \G$ and $f$ may be any object upon which adjunction makes sense.  For examples, $f \in \G$ or $f \subseteq \G$ or $f$ is a representation into $\G$.

The commutator map is
\begin{equation*}
\CC : \G \times \G \to \G, \ \ \ \ \ \CC(g,h) = ghg^{-1}h^{-1}.
\end{equation*}
\end{defin}
For $1 \le i \le 2p$, let $A_i \in \pi$ be a simple non-separating loop such that the intersection number of $A_i$ and $A_j$ equals the Kronecker $\delta_{(i)(j-p)}$.  Then there is a presentation (see Figure 1)
\begin{equation}\label{eqn:standard presentation}
\pi  =  \langle A_i, \ 1 \le i \le 2p \ | \ \prod_{i=1}^{p} \CC(A_i,A_{i+p}) \rangle.
\end{equation}

%

Suppose $\G$ is a group.  Denote by $Z(\G)$
and $\Aut(\G)$ the center
and the automorphism group of $\G$, respectively.  If $\G$ acts on a set $S$, denote by $S^\G$ the fixed point subset of $S$.  If $a \in \G$, then $S^a$ means $S^{\langle a \rangle}$, where $\langle a \rangle$ is the cyclic group generated by $a$.

%

For the rest of the paper, let $G$ be a complex reductive group and $\I$ the identity element of $G$.  Let $d, r$ and $z$ denote the dimension and rank of $G$ and the dimension of $Z(G)$, respectively.  Let
$H$ and
$K$ be a Cartan and a maximum compact subgroup of $G$ with Lie algebras $\fh,
\fk$ and $\fg$, respectively.
Let $N$ be the set of torsion points of $G$.
\begin{defin}\label{defin:irreducible}
Suppose $\rho : \G \to G$ is a group homomorphism.  Then $\rho$ is {\bf irreducible} if $\rho(\G)$ is not contained in any proper parabolic subgroup of $G$.
\end{defin}

Write $R(G)$ (resp. $R_0(G)$) for the representation variety $\Hom(\pi,G)$ (resp. $\Hom(\pi_1(X_0),G)$).  Then
\begin{equation}\label{eq:variety}
R(G)  =  \{a \in G^{2p} : \prod_{i=1}^{p} \CC(a_i,a_{i+p}) = \I\}.
\end{equation}
In other words, $R(G)$ may be identified with a subvariety of $G^{2p}$.  The group $G$ acts on $R(G)$ by conjugation
$$
G \times R(G) \longrightarrow R(G), \ \ (g,\rho) = g.\rho = g \rho g^{-1}.
$$
The resulting categorical quotient is $R(G) \longrightarrow \HH^1(X,G)$.

\subsection{The mapping class groups and their actions}\label{sec:MCG}

Let $Diff(X)$ be the group of orientation preserving diffeomorphisms of $X$ and $Diff(X,x)$ the subgroup fixing the base point $x$.  There is a natural inclusion $\iota : Diff(X,x) \hookrightarrow Diff(X)$.  Let $\Gamma = \pi_0(Diff(X))$ and $\Gamma(x) = \pi_0(Diff(X,x))$.  Since every component of $Diff(X)$ contains an element that fixes $x$, the induced map $\iota_* : \Gamma(x) \to \Gamma$ is onto.
The group $\Gamma(x)$ acts on $\pi(x)$:
$$
\Gamma(x) \times \pi(x) \longrightarrow \pi(x), \ \ \ \
(\gamma, u) \longmapsto \gamma(u).
$$
This induces a $\Gamma(x)$-action on
$R(G)$:
\begin{equation} \label{Gamma2}
\Gamma(x) \times R(G) \longrightarrow
R(G), \ \ \ \ (\gamma, \rho) \longmapsto \rho
\circ \gamma^{-1}.
\end{equation}
This action factors through $\iota_*$ and induces a $\Gamma$-action on $\HH^1(X,G)$.

\begin{defin}
Suppose $L$ is a simple loop on $X$ such that $x \not\in L$.  Let $\B$ be a small tubular neighborhood of $L$ with $x \not\in \B$.  Cutting along $L$, there are two ways to Dehn twist along $L$ in $\B$.  Stay consistently on one side of $X$ and call the right turn direction {\bf positive} \cite{ES1} and denote by $\tau_{L}$ the Dehn twist in the {\bf negative} direction along $L$.
\end{defin}
The construction of $\tau_{L} \in Diff(X,x)$ actually depends on $\B$ and the details of the twist; however, its images in $\Gamma(x)$ and $\Gamma$ under $\iota_*$ do not.   Denote also by $\tau_{L}$ its images in $\Gamma(x)$ and $\Gamma$.

\begin{defin}
Suppose $C \in \pi$ is simple.  Let $C'$ denote a simple loop that is a deformation of $C$ in a small tubular neighborhood of $C$ such that $x \not\in C'$.

In the cases of $A_i \in \pi$ in the standard presentation \eqref{eqn:standard presentation}, deform $A_i$ slightly left and call the resulting loop $A_i'$ for each $1 \le i \le 2p$ as in Figure 1.
\end{defin}

By a direct calculation \cite{G1, G2} (see Figure 1),
\begin{prop}\label{Dehn twist}
Let $\rho \in R(G)$.
\begin{enumerate}
\item
Let
$C \in \pi$ be simple and separating and $X_1, X_2$ the resulting connected components of $X \setminus C'$.  Then the Dehn twist $\tau_{C'}$-action on $\HH^1(X,G)$ is covered by the $\tau_{C'}$-action on $R(G)$ with
    $$
    \tau_{C'}(\rho)(A) = \left\{
    \begin{array}{lll}
    \rho(A) & \text{ if } & A \in \pi_1(X_1)\\
    \rho(B A B^{-1}) & \text{ if } & A \in \pi_1(X_2)
    \end{array}\right.,
    $$
   where $B = C$ or $B = C^{-1}$.
\item The $\tau_{A_i'}$-action and the $\tau_{A_{p+i}'}$-action on $\HH^1(X,G)$ are covered by
$$
\tau_{A_i'}(\rho)(A_{p+i})  =  \rho(A_{p+i} A_i^{-1}), \ \ \
\tau_{A_{p+i}'}(\rho)(A_i) =  \rho(A_i A_{p+i}),
$$
for $1 \le i \le p$.
\end{enumerate}
\end{prop}

\section{Dehn twists and degenerations}
Consider a degeneration $\f : \X \to D$.  The map $c^*$ is defined as follows.  The map
$$
c : X \to X_0
$$
induces a homomorphism on the fundamental group
$$
c_* : \pi \to \pi_1(X_0).
$$
This further induces a map
$$
c^* : R_0(G) \to R(G), \ \ \ c^*(\rho) = \rho \circ c_*
$$
which descends to the map $c^*$ of \eqref{cT}.
Hence
\begin{lemma}\label{lem:pinched cycle}
If $[\rho] \in \im(c^*)$, then $\rho(C) = \I$ for all $C \in \ker(c_*)$.
\end{lemma}

\begin{defin}\label{def:pseudo-degeneration}
Let $\S = \{C_1, \cdots, C_m\} \subseteq \pi$ be a set of simple loops.  Suppose there exists a set of pairwise disjoint simple loops $\{C_1', \cdots, C_m'\}$ such that $C_i'$ is a deformation of $C_i$ in a small tubular neighborhood of $C_i$.
Then the element $$\tau_\S = \prod_{i = 1}^m \tau_{C_i'}^{n_i} \in \Gamma$$ is called a {\bf pseudo-degeneration} for $\S$.
We always assume that $n_i \in \Z \setminus \{0\}.$
An element $[\rho] \in \HH^1(X,G)$ is {\bf exceptional} for $\tau_\S$ if $\tau_\S([\rho]) = [\rho]$ and $\rho(C_i) \neq \I$ for some $i$.
\end{defin}
The construction in \cite[Section 2]{ES1} and \cite[Theorem 3]{ES1} completely determines when a pseudo-degeneration arises from a degeneration:
\begin{thm}[Earle, Sipe]\label{thm:degenerate}
Let $\tau_\S = \prod_{i=1}^m \tau_{C_i'}^{n_i} \in \Gamma$ be a pseudo-degeneration for $\S = \{C_1, \cdots, C_m\} \subseteq \pi$.  Then $\tau_\S = T_*$ for some degeneration $\f : \X \to D$ with $\S \subseteq \ker(c_*)$ if and only if $n_i > 0$ for all $i$.
\end{thm}

\begin{rem}\label{rem:local degeneration}
A degeneration as defined in \cite{ES1} allows the central fibre to have general isolated singularities.  However the actual construction in Section 2 of \cite{ES1} clearly shows that the central fibre has normal crossings and the degeneration is of the type $\{u v = t^{n_i}\} \to t$ near the singularity associated with $C_i$.
\end{rem}

By Lemma \ref{lem:pinched cycle},
\begin{cor}\label{cor:DehnTwist}
Suppose that $\f : \X \to D$ is a degeneration as in Theorem \ref{thm:degenerate} and $\tau_\S = T_*$ for $\S = \{C_1, \cdots, C_n\}$.  If $[\rho] \in \im(c^*)$, then $\rho(C_i) = \I$ for all $i$.
\end{cor}

\section{Exceptional families}\label{sec:exceptional}

In this section, we prove that the exceptional loci $\V_\f$ can be  rather large. Recall that $N \neq \{\I\}$ is the set of torsion points of $G$.

\begin{thm} \label{thm:exception}
\mbox{ }
\newline
\begin{enumerate}
\item Suppose $p = 1$.  Then there exists a degeneration $\f : \X \to D$ with $\dim(\V_\f) \ge r$.

\item Suppose $p > 1$. Then there exists a degeneration $\f : \X \to D$ and $[\rho] \in \V_\f$ with $\im(\rho)$ being irreducible and Zariski dense in $G$.  Moreover $\dim(\V_\f) \ge (2p-3)d + 2z$.
\end{enumerate}
\end{thm}
\begin{proof}
Suppose that $\lambda \in N \setminus \{\I\}$ is an $n$-torsion point and
$$
V_\lambda = \{\rho \in R(G) : \rho(A_1) = \lambda\} = \{a \in R(G) : a_1 = \lambda\},
$$
via the identification \eqref{eq:variety}.  Let $\rho \in V_\lambda$. By Theorem \ref{thm:degenerate}, there is a degeneration $\f : \X \to D$ with $T_* = \tau_\S = \tau_{A_1'}^n$, where $\S = \{A_1\} \subseteq \ker(c_*)$.  Since $\rho(A_1)^{-n} = \lambda^{-n} = \I$, $\tau_\S(\rho) = \rho$ by Proposition \ref{Dehn twist}(2).  Hence $[\rho]$ is a cocycle class invariant under $T^*$.  Since $\rho(A_1) = \lambda \neq \I$, $[\rho] \not\in \im(c^*)$ by Corollary \ref{cor:DehnTwist}.  Hence $\pr(V_\lambda) \subseteq \V_\f$.

Suppose $p=1$.  Then
$$V_\lambda \supseteq \{a_2 \in G : \CC(\lambda,a_2) = \I \} = G^\lambda.$$
Then $\dim(\V_\f) \ge r.$ This proves part (1).

Suppose $p > 1$.
There exists a two generator subgroup $\langle g ,h\rangle$ that is irreducible and Zariski dense in $G$.  Since $N \cap K$ is dense in $K$ in the usual topology, we may choose $\lambda \in N \cap K$ and $h \in G$ such that $\langle \lambda ,h\rangle$ is irreducible and Zariski dense in $G$.  Since $p > 1$, there exists $\rho \in V_\lambda$ such that
$$
\left\{
\begin{array}{l}
\rho(A_1) = \rho(A_{p+2}) = \lambda, \ \  \rho(A_2) = \rho(A_{p+1}) = h\\
\rho(A_i) = \rho(A_{p+i}) = {\I}, \ \   \text{ for } 3 \le i \le p
\end{array}
\right..
$$
Then $\rho$ satisfies the required condition.
A direct calculation then shows that
$$
\dim(V_\lambda) \ge (2p - 2)d + z, \ \ \
\dim(\V_\f) \ge (2p-3)d + 2z.
$$
\end{proof}

\section{More simple examples}

The previous section shows that the exceptional family can be rather large even when $G$ is abelian.  In this section, we give more examples corresponding to a more restrictive class of degenerations.  For simplicity, assume that $G$ is semi-simple.

\begin{defin}\label{def: simple degeneration}
A degeneration $\f : X \to D$ is {\bf simple} if it gives rise to a pseudo-degeneration with $n_i = 1$ for all $i$, where $n_i$ is as in Theorem \ref{thm:degenerate}.
\end{defin}

\begin{lemma} \label{lem:Csurj}
The map $\CC : K \times K \to K$ is surjective; hence,
the center $Z(G)$ is contained in $\CC(K \times K)$.
\end{lemma}
\begin{proof}
Since $G$ as well as $K$ are semi-simple, $\CC(K \times K) = K$ by \cite[Remark (2.1.5)]{PX1}.  The Lemma follows because $Z(G)$ is finite and contained in $Z(K)$.
\end{proof}

The following Proposition \ref{prop:zariski open} is a direct consequence of \cite[Theorem 1.3]{La1} as pointed out to me by Jiu-Kang Yu:
\begin{prop}\label{prop:zariski open}
If $p \ge 1$, then the space
$$
\{(a_1,\cdots,a_{2p}) \in G^{2p} : \langle a_1, \cdots, a_{2p}\rangle \text{ is Zariski dense in } G \}
$$
is Zariski open in $G^{2p}$.
\end{prop}
By identification \eqref{eq:variety}, $R(G)$ is a subvariety of $G^{2p}$.  Since the intersection of a subvariety with a Zariski open set is Zariski open in the subvariety, we have
\begin{cor}\label{cor:zariski open}
If $p \ge 1$, then the space
$$
\{\rho \in R(G) : \im(\rho) \text{ is Zariski dense in } G \}
$$
is Zariski open in $R(G)$.
\end{cor}

\begin{thm} \label{thm:simple exception}
\mbox{ }
\newline
\begin{enumerate}
\item If $p > 0$, then there exists simple degeneration $\f : \X \to D$ such that $\V_\f \neq \emptyset$.

\item If $p > 1$ and $G$ contains a closed semi-simple subgroup $G'$ with $Z(G') \neq \{\I\}$ and $\dim(G') = d'$, then there exists simple degeneration $\f : \X \to D$ such that
$\dim(\V_\f) \ge (2p-2)d' - d$.

\item If $p > 2$ and $G$ contains a closed semi-simple subgroup $G'$ with $Z(G') \neq \{\I\}$ and $\dim(G') = d'$, then there exists simple degeneration $\f : \X \to D$ with $[\rho] \in \V_\f$ such that $\im(\rho)$ is irreducible and Zariski dense in $G'$;  moreover, $\dim(\V_\f) \ge (2p-3)d'$.
\end{enumerate}
\end{thm}
\begin{proof}
Let $1 \le n \le p$ and let $N(H)$ be the normalizer of $H$.
Consider the adjoint action of $N(H)$ on $H$
$$
N(H) \times H \to H, \ \ (w,h) = w. h = whw^{-1}.
$$
Suppose that there exist $w \in N(H)$ and $g, h \in H$
such that $g \neq \I$ and
\begin{equation}\label{eq:weyl equation}
\left\{
\begin{array}{lll}
w.(hg) & = & h\\
\CC(g,w)^n & = & \I
\end{array}\right..
\end{equation}
Then
\begin{equation}\label{eq:conjpair}
h.(g,wg^{-1}) = (hgh^{-1}, hwg^{-1}h^{-1}) = (g,w).
\end{equation}
Let $1 \le i_1 < \cdots < i_n \le p$ and
$\rho \in R(G)$ such that
\begin{equation}\label{eq:def}
    \rho(A_l) = \left\{
    \begin{array}{lll}
    g & \text{ if } & l = i_j \text{ for some } j\\
    w & \text{ if } & l = i_j + p \text{ for some } j\\
    \I & \text{ if } & \text{ otherwise }
    \end{array}\right.
\end{equation}
By Proposition \ref{Dehn twist}(2) and Theorem \ref{thm:degenerate}, there exists a degeneration $\f : \X \to D$ such that $T_* = \tau_\S = \prod_{j=1}^n \tau_{A_{i_j}'}$ with $\S = \{A_{i_1}, \cdots, A_{i_n}\} \subseteq \ker(c_*)$.  By \eqref{eq:conjpair} and Proposition \ref{Dehn twist}(2), $\tau_\S(\rho) = h.\rho$.  Hence, $\tau_\S([\rho]) = [\rho]$, i.e. $[\rho] \in \HH^1(X,G)^{T^*}$.  Since $\rho(A_{i_j}) = g \neq \I,$ $[\rho] \not\in \im(c^*)$ by Corollary \ref{cor:DehnTwist}.  Hence $[\rho] \in \V_\f$.

For any $g \in H$, since $w.g \in H$, $w.g$ commutes with $g$.  Hence $\CC(g,w)^n = g^n(w.g)^n$.  Hence $\CC(g,w)^n = \I$ if and only $(w.g)^n = g^{-n}$.

Since $G$ either contains $\PGL(2)$ or $\SL(2)$, part (1) will follow if we find triples $w, g, h$ satisfying Equation \eqref{eq:weyl equation} with $g \neq \I$ for the cases of $G = \PGL(2)$ and $G = \SL(2)$.

Let $1 \le n \le p$ and $G = \PGL(2)$ (resp. $G = \SL(2)$).  For $1 \le i_1 < \cdots < i_n \le p$, let $s \in \C \setminus \{\pm 1\}$ (resp. $s \in \C \setminus \{1\}$) with $s^{2n} = 1$ (resp. $s^n = 1$) and
$$
g = \bmatrix s & 0 \\ 0 & 1/s  \endbmatrix, \ \
w = \bmatrix 0 & 1 \\ -1 & 0  \endbmatrix, \ \
h = \bmatrix t & 0 \\ 0 & 1/t  \endbmatrix,
$$
where $t^2 = s$.
Let $\rho \in R(G)$ be as defined by \eqref{eq:def}.
Then $h.(g,wg^{-1}) = (g, w)$.  Hence $[\rho] \in \V_\f$ with  $A_{i_j} \in \ker(c_*)$.
Incidentally, $\im(\rho)$ is abelian in these two cases if $p = 1$, but observe that $h \not\in \im(\rho)$.
This proves part (1).
Notice that in this case $\tau_\S(\rho) \neq \rho$.

For part (2) and (3), let $K'$ be the maximum compact subgroup of $G'$.

Suppose $p > 1$.
Let $\lambda \in Z(G') \setminus \{\I\}$ and
$$
V_\lambda = \{\rho \in R(G') : \rho(\CC(A_1,A_{p+1})) = \lambda\} = \{a \in R(G') : \CC(a_1,a_{p+1}) = \lambda\},
$$
via the identification \eqref{eq:variety}.
Since $p > 1$, by Lemma \ref{lem:Csurj}, $V_\lambda \neq \emptyset$.  A direct calculation shows that
$$
\dim(V_\lambda) \ge  2(p-1)d'.
$$
Let $\rho \in V_\lambda$.  The element $\CC(A_1,A_{p+1})$ corresponds to a simple and separating loop $C$ on $X$.  By Theorem \ref{thm:degenerate}, there is a degeneration $\f : \X \to D$ such that $T_* = \tau_\S = \tau_{C'}$ with $\S = \{C\}$.  Since $\rho(C) = \lambda \in Z(G')$, $\tau_\S(\rho) = \rho$ by Proposition \ref{Dehn twist}(1).  Hence $[\rho]$ is a cocycle class invariant under $T^*$.  Since $\rho(C) = \lambda \neq \I$, $[\rho] \not\in \im(c^*)$ by Corollary \ref{cor:DehnTwist}.  Hence $\pr(V_\lambda) \subseteq \V_\f$.
Then
$$
\dim(\V_\f) \ge \dim(V_\lambda) - d \ge 2(p-1)d' - d.
$$
This proves part (2).

Suppose in addition that $p > 2$.
Let $g,h \in K'$ such that $\langle g, h\rangle$ is Zariski dense in $G'$.  Then $(\lambda\CC(g,h))^{-1} \in \CC(K' \times K')$ by Lemma \ref{lem:Csurj}.  Hence there exits $\rho_0 \in V_\lambda$  such that
\begin{equation*}
\left\{
\begin{array}{l}
\rho_0(\CC(A_1,A_{p+1}))  =  \lambda,\\
\rho_0(A_2) = g, \ \ \rho(A_{p+2}) = h, \\
\rho_0(\CC(A_3,A_{p+3})) = \lambda^{-1} \CC(g,h)
\end{array}\right.
\end{equation*}
Then $\langle g,h \rangle \subseteq \im(\rho) \subseteq G'$.   Hence $\pr(V_\lambda) \subseteq \V_\f$.

The property of $\im(\rho)$ being irreducible is open.
Since $G'$ is semi-simple, by Corollary \ref{cor:zariski open},
there exists a Zariski open set $V_0 \subseteq V_\lambda$ containing $\rho_0$ such that $\rho \in V_0$ if and only if $\rho$ is irreducible and with $\im(\rho)$ being Zariski dense in $G'$.

Let $\rho \in V_0$ and
$$G_\rho = \{g \in G:  g.\im(\rho) \subseteq G'\}.$$
Let $g \in G_\rho$.
Since $\im(\rho)$ is Zariski dense in $G'$, $g.G' = G'$, i.e. the $G_\rho$-action on $R(G')$ factors through $\Aut(G')$.  Since $G'$ is semi-simple, $\dim(\Aut(G')) = d'$.  Hence the $G_\rho$-orbit of $\rho$ is at most $d'$-dimensional.  Hence
$$
\dim(\V_\f) \ge \dim(V_0) - d' \ge (2p - 3)d'.
$$
Part (3) now follows.

\end{proof}

\begin{rem}\label{rem:commutatorimage}
For an abelian $G$, the $\Ad(G)$-action on itself is trivial.  However, when $G$ is semi-simple, the $N(H)$-action on $H$ is not trivial and allows the construction of the exceptional families in Theorem \ref{thm:simple exception}(1).
When $G$ is abelian, $\CC(G \times G) = \{\I\}$ while almost the opposite is true when $G$ is semi-simple in the sense that $\CC(G \times G)$ is Zariski dense and contains $K$.  This allows one to construct the exceptional families in Theorem \ref{thm:simple exception}(2, 3).
\end{rem}

The examples we have constructed motivate the following conjecture on local invariant classes:
\begin{conjecture}\label{conj:lic}
Suppose $\f : \X \to D$ is a degeneration.
If
$
[\rho] \in \HH^1(X,G)^{T^*},
$
then $\rho(C) \in N$ for all simple $C \in \ker(c_*)$.
\end{conjecture}
\begin{rem}
Notice that even if $G$ is reductive and not necessarily semi-simple,  it is still the case that $\CC(G,G) \cap Z(G) \subseteq N$; otherwise, the construction for Theorem \ref{thm:simple exception}(2) would have yielded counterexamples to Conjecture \ref{conj:lic}.
\end{rem}

\section{Pseudo-degenerations}
We do not have a counter-example for Conjecture \ref{conj:lic} (obviously).  However we do have the following suggestive example.  Assume $G$ is semi-simple in this section.

%
%
%

\begin{lemma}\label{lem:onto}
Let $g \in K$.  Then
$$
\CC : K \times K^g \to K
$$
is onto.
\end{lemma}
\begin{proof}
It is sufficient to assume $g$ is regular.  The Lemma then follows from \cite[Proposition (B.1)]{PX1}.
\end{proof}

%


Let ${\rm pr}_1 : W' \to G$
be the projection to the first factor, where
$$
W' = \{(g,h,k) \in G \times G \times G: k.(g,hg^{-1}) = (g,h)\}.
$$
\begin{prop}\label{prop:dim}
The image ${\rm pr}_1(W')$ contains $K$.
\end{prop}
\begin{proof}
The equation $k.(hg^{-1}) = h$ is equivalent to $g = \CC(k^{-1},h^{-1}).$
Let $g \in K$.
By Lemma \ref{lem:onto}, $${\rm pr}_1^{-1}(g)  = \{(h,k) \in G \times G^g : g = \CC(k^{-1},h^{-1})\} \neq \emptyset.$$

\end{proof}

Consider the presentation of $\pi$ as in Figure 1.  Let $W \subseteq R(G)$ such that $\rho \in W$ if
\begin{equation*}
\left\{
\begin{array}{l}
\rho(A_1) = \rho(A_{p+2}) = g\\
\rho(A_2) = \rho(A_{p+1}) = h\\
\rho(A_i), \rho(A_{p+i}) \in G^k \text{ for } 3 \le i \le p
\end{array}\right.,
\end{equation*}
where $(g,h,k) \in W'$.
Let
$$\S = \{A_1, A_{p+2}\}, \ \ \ \tau_\S = \tau_{A_1'} \tau_{A_{p+2}'}^{-1}.$$
Let $\rho \in W$.  Then by Proposition \ref{Dehn twist}(2) (see Figure 1),
\begin{eqnarray*}
\tau_\S(\rho)(A_2) = \rho(A_2 A_{p+2}^{-1}) = hg^{-1} = \rho(A_{p+1} A_1^{-1}) = \tau_\S(\rho)(A_{p+1});
\end{eqnarray*}
moreover, $\tau_\S(\rho)(A_i) = \rho(A_i)$ for $1 \le i \le 2p$ and $i \neq 2, p+1$.
By Proposition \ref{prop:dim},
$\tau_\S(\rho) = k.\rho$; hence, $\tau_\S([\rho]) = [\rho]$.  Moreover,
the image of the map
$$
W \to G : \ \ \ \rho \to \rho(A_1)
$$
contains $K$.  Hence the family $\pr(W)$ would have contained counter-examples to Conjecture \ref{conj:lic} if $\tau_\S$ had actually arisen from a degeneration.

For an explicit construction with $G = \SL(2)$, let $s, t \in \C^\times, t^2 = s$ and
$$
g = \bmatrix s & 0 \\ 0 & 1/s  \endbmatrix, \ \
h = \bmatrix 0 & 1 \\ -1 & 0  \endbmatrix, \ \
k = \bmatrix t & 0 \\ 0 & 1/t  \endbmatrix.
$$

Since $\tau_\S$ contains a negative power $\tau_{A_{p+1}}^{-1}$, by Theorem \ref{thm:degenerate}, the pseudo-degeneration $\tau_\S$ does not correspond to a degeneration.  However it can be realized by the following family of hyperelliptic curves parameterized by $t$:
$$
y^2 = (z + (t-1)a_2 - ta_1) (z - a_2)(z - a_3)(z + ({\bar t}-1)a_5 - {\bar t} a_4)\prod_{i=5}^{2g+2} (z - a_i),
$$
where $a_1, \cdots, a_{2p+2}$ are distinct complex numbers and ${\bar t}$ means the complex conjugate of $t$.

\section{The case of $\tSL$}\label{sec:tSL}
In this final section, we give yet another suggestive example with regard to Conjecture~\ref{conj:lic}.  Up to this point, $G$ is assumed to be complex algebraic and this is a standing assumption for Conjecture~\ref{conj:lic}.  Now we consider $\tSL$, the universal cover of $\SL(2,\R)$.  The group $\tSL$ is not algebraic; however, one may still define $R(\tSL)$ as a topological space (see Definition (\ref{eq:variety})) and $\HH^1(X,\tSL)$ as its topological quotient by the adjoint $\tSL$-action.  Similarly, there is a $\Gamma$-action on $R(\tSL)$ which descends to a $\Gamma$-action on $\HH^1(X,\tSL)$.  Likewise, if $\f : \X \to D$ is a degeneration, then the corresponding map $T^* : \HH^1(X,\tSL) \to \HH^1(X,\tSL)$ is similarly defined.

The center $Z(\tSL)$ is isomorphic to $\Z$.  Let $\eta \in Z(\tSL)$ be a generator and $\CC_2 = \{ab : a,b \in \CC(\tSL \times \tSL)\}.$  Then $\eta, \eta^{-1} \in \CC_2$ (see \cite{Go3}, Sections 9 \& 10).
\begin{prop}
If $p > 3$, then there exists a degeneration $\f : \X \to D$, a simple $C \in \ker(c_*)$ and $[\rho] \in \HH^1(X,\tSL)^{T^*}$ such that $\rho(C)$ has infinite order.
\end{prop}
\begin{proof}
The proof is similar to the proof of Theorem~\ref{thm:simple exception} (2).  Let $C = \prod_{i=1}^2 \CC(A_i,A_{p+i}).$  Then $C$ is represented by a simple and separating loop on $X$.
Since $\eta, \eta^{-1} \in \CC_2$ and $p > 3$, there exits $\rho \in R(\tSL)$ such that
$$\rho(C) = \rho(\prod_{i=1}^2 \CC(A_i,A_{p+i})) = \eta, \ \ \rho(\prod_{i = 3}^p \CC(A_i,A_{p+i})) = \eta^{-1}.$$
Therefore the order of $\rho(C)$ is infinite.
By Theorem \ref{thm:degenerate}, there is a degeneration $\f : \X \to D$ such that $T_* = \tau_\S$ with $\S = \{C\}$.  Since $\rho(C) = \eta \in Z(\tSL)$, by Proposition~\ref{Dehn twist} (1), $\tau_\S(\rho) = \rho$ which means that $[\rho] \in \HH^1(X,\tSL)^{T^*}$.
\end{proof}




\begin{thebibliography}{99}

\bibitem{A1}
Andersen, J{\o}rgen Ellegaard, Fixed points of the mapping class group in the ${\rm SU}(n)$ moduli spaces. {\em Proc. Amer. Math. Soc.} {\bf 125} (1997), no. 5, 1511--1515.


\bibitem{C1}
 Clemens, C. H., Degeneration of K\"ahler manifolds. {\em Duke Math. J.} {\bf 44} (1977), no. 2, 215--290.



\bibitem{ES1}
Earle, Clifford J.; Sipe, Patricia L., Families of Riemann surfaces over the punctured disk.  {\em Pacific J. Math.}  {\bf 150}  (1991),  no. 1, 79--96.


\bibitem{G1} Goldman, William M., Invariant functions on Lie groups and Hamiltonian flows of surface group representations.  {\em Invent. Math.}  {\bf 85}  (1986),  no. 2, 263--302.

\bibitem{Go3} Goldman, William M., Topological components of spaces of representations. {\em Invent. Math.} {\bf 93} (1988), no. 3, 557-–607.

\bibitem{G2}
Goldman, W. M., Ergodic theory on moduli spaces.  {\em Ann. of Math.} (2)  {\bf 146}  (1997),  no. 3, 475--507.




\bibitem{La1}
Larsen, Michael, Rigidity in the invariant theory of compact groups. {\em arXiv:math/0212193}.

\bibitem{LM1}
Lubotzky, A.; Magid, A. R., Varieties of representations of finitely generated groups.  {\em Mem. Amer. Math. Soc.}  {\bf 58}  (1985),  no. 336, xi+117 pp.

\bibitem{M1}
Mumford, D.; Fogarty, J.; Kirwan, F. Geometric invariant theory. Third edition. Ergebnisse der Mathematik und ihrer Grenzgebiete (2) [Results in Mathematics and Related Areas (2)], {\bf 34}. {\em Springer-Verlag, Berlin,} 1994. xiv+292 pp. ISBN: 3-540-56963-4

\bibitem{PS1}
Peters, Chris A. M.; Steenbrink, Joseph H. M., Mixed Hodge structures. Ergebnisse der Mathematik und ihrer Grenzgebiete, {\bf 52} {\em Springer-Verlag, Berlin}, 2008. xiv+470 pp.


\bibitem{PX1}
Pickrell, Doug; Xia, Eugene Z., Ergodicity of mapping class group actions on representation varieties. I. Closed surfaces. {\em Comment. Math. Helv.} {\bf 77} (2002), no. 2, 339--362.



\bibitem{S1}
Schmid, Wilfried, Variation of Hodge structure: the singularities of the period mapping. {\em Invent. Math.} {\bf 22} (1973), 211--319.

\bibitem{S2}
Steenbrink, Joseph, Limits of Hodge structures. {\em Invent. Math.} {\bf 31} (1975/76), no. 3, 229--257.

\bibitem{TX1}
Tsai, Yen-Lung; Xia, Eugene Z., Non-abelian local invariant cycles. {\em Proc. Amer. Math. Soc.} {\bf 135} (2007), no. 8, 2365--2367.



\end{thebibliography}
\end{document}